\input amstex
\input epsfx.tex
\frenchspacing
\documentstyle{amsppt}
\magnification=\magstep1
\baselineskip=14pt
\vsize=18.5cm
\footline{\hfill\sevenrm version 20131201} 
 
\def\Tr{{\text{\rm Tr}}}
\def\Gal{{\text{\rm Gal}}}

\def\Hom{{\text{\rm Hom}}}

\def\SL{{\text{\rm SL}}}
\def\GL{{\text{\rm GL}}}

\def\Ker{\mathop{\text{\rm ker}}}

\def\legendre{\overwithdelims()}
\def\endproof{$\square$}

\def\mapright#1{\ \smash{\mathop{\longrightarrow}\limits^{#1}}\ }

\def\O{{\Cal O}}
\def\Q{\bold Q}
\def\C{\bold C}
\def\R{\bold R}

\def\Z{\bold Z}

\def\gp{{\goth p}}

\newcount\refCount
\def\newref#1 {\advance\refCount by 1
\expandafter\edef\csname#1\endcsname{\the\refCount}}
\newref CF 
\newref Var 
\newref Coh 
\newref EDiss 
\newref EP 
\newref Hoff 
\newref KP 
\newref Ku 
\newref Neu 
\newref PatCub 
\newref PatIJM 
\newref PatDis 
\newref Suz 
\newref Tate 
\newref Well 
\topmatter
\title 
Fourier coefficients of sextic theta series
\endtitle
\author Reinier Br\"oker, Jeff Hoffstein\endauthor
\address Brown University, Department of Mathematics,
Box 1917, Providence, RI, USA
\email reinier, jhoff\@math.brown.edu\endemail
\endaddress
\abstract
This article focuses on the theta series on the 6-fold cover of $\GL_2$. We 
investigate the Fourier coefficients $\tau(r)$ of the theta series, and give
partially proven, partially conjectured values for $\tau(\pi)^2$, 
$\tau(\pi^2)$ and $\tau(\pi^4)$ for prime values $\pi$.
\endabstract
\endtopmatter

\document

\head 1. Introduction
\endhead
\noindent
The Jacobi theta function, defined by
$$
\theta(z) = \sum_{n= -\infty}^\infty e^{2 \pi i n^2 z},
$$
for $z = x+iy$ with $y>0$ is of fundamental importance in many areas of 
mathematics. In many applications, $x$ is set equal to 0, $y$ is set 
equal $t/2$ and the relevant property of the theta function is the 
transformation formula
$$
 \sum_{n= -\infty}^\infty e^{- \pi n^2 /t} = \sqrt{t} \sum_{n= -\infty}^\infty e^{- \pi n^2 t}.
$$
However, $\theta(z) $ possesses a more general transformation property. We let 
$$
\Gamma_0(4) = \{ \gamma\in \SL_2(\Z) \mid \gamma = \binom{a\thinspace\, b}{c\thinspace \,d}, \text{with}\,\, c \equiv 0 \bmod 4 \}
$$
be the usual congruence subgroup.
Then, for $\gamma \in \Gamma_0(4)$ we have
$$
\theta (\gamma z) = j(\gamma,z)\theta(z),
$$
where 
$$
 j(\gamma,z) = \epsilon_d^{-1}\left(\frac c d\right)\sqrt{cz+d}.
$$
Here 
$\epsilon_d =  1$ if $d \equiv 1 \bmod 4$, and $\epsilon_d =  i$ 
if $d \equiv 3 \bmod 4$, and $\left(\frac c d\right)$ is the usual 
quadratic symbol except that we multiply by $-1$ for $c,d<0$.
The square root $\sqrt{cz+d}$ is chosen to 
have argument with absolute value less than $\pi/2$.

The theta function has a beautiful connection with Eisenstein series of 
half-integral weight. One can construct such an Eisenstein series as follows:
$$
E^{(2)}(z,s) = \sum_{\Gamma_\infty \backslash \Gamma_0(4)}\text{Im}(\gamma z)^s \frac{\theta(z)}{\theta(\gamma z)}
= \sum_{(c,d)=1, c \ge 0\atop  c \equiv 0 \bmod 4}\frac{\epsilon_d \left(\frac c d\right)y^{s}}{|cz +d|^{2s}\sqrt{cz+d}}.
$$
This converges absolutely for $\text{Re}(s) >3/4 $ and, by construction, 
satisfies the same transformation property as $\theta(z)$, namely
$$
E^{(2)}(\gamma z,s) = j(\gamma,z)E^{(2)}( z,s).
$$
The remarkable thing is that $E^{(2)}( z,s)$ has a simple pole at $s = 3/4$, 
and upon taking the residue, one recovers the original theta function.  In 
other words, the equality
$$
\text{\rm Res}_{s = 3/4} E^{(2)}( z,s) = c \theta(z),
$$
holds for some non-zero constant $c$.

Weil made the observation that just as an automorphic form on the upper half 
plane can be interpreted as an automorphic form on the group $G= GL_2(\R)$, 
the functions $\theta(z)$ and $E^{(2)}( z,s)$ can be interpreted as functions 
on $\tilde G$, the 2-fold metaplectic cover of $G$.   Here, one has 
$$
\tilde G = \{ (g,\epsilon)| g \in G, \epsilon = \pm 1\},
$$
and multiplication is defined by
$$
(g,\epsilon)(g',\epsilon') = (gg', \epsilon \epsilon' \sigma(g,g')),
$$
with $\sigma(g,g')$ a certain explicit 2-cocycle.

Kubota, [\Ku], defined Eisenstein series on the $n$-fold metaplectic cover 
of $\GL_2$, and observed that these Eisenstein series have simple poles 
at $s = 1/2 + 1/(2n)$. The residues at this point are automorphic forms on 
the $n$-cover of $\GL_2$, and generalize the notion of the quadratic theta 
function. Unlike the quadratic theta function however, the Fourier 
coefficients of the generalized theta function when $n \ge 3$ are very 
mysterious, and at present are only completely understood in the case $n=3$.
In this introduction we will survey what is known and conjectured about these 
theta functions.   To make the underlying structure clearer we will  be very 
imprecise in this in the sense that we will ignore bad primes and assume that 
reciprocity works perfectly.

Kubota's Eisenstein series can be defined in the following way in the 
cases $n=3,4,6$. We let $\zeta_n$ be a primitive $n$-th root of unity, and put $K=\Q(\zeta_n)$ with ring of integers $\Z[\zeta_n]$.
Let $\left(\frac c d\right)_n$ represent the $n$-th order residue symbol 
and, for 
$$
z = \binom{y\thinspace\, x}{0\thinspace \,1}k \binom{\alpha\thinspace\, 0}{0\thinspace \,\alpha}
$$
with $x \in \C$ and $y>0$, and $k \in U(2,\C)$, $\alpha  \in \C^*$, 
let $I(z) = y$. Kubota observed that for suitable $N$, the $n^{th}$ power 
reciprocity law implies that the function
$$ 
\kappa(\gamma) = \left(\frac c d\right)_n,
$$
from $\Gamma(N) =  \Gamma_N(\SL_2(\Z[\zeta_n]))\rightarrow \C^*$, is a homomorphism. (The
choice $N=n^2$ works, but need not be minimal.)
   He used this to define the Eisenstein series
$$
E^{(n)}(z,s) = \sum_{\Gamma_\infty \backslash \Gamma(N)}\kappa(\gamma)I(\gamma z)^{2s},
$$
which converges absolutely for $\text{Re}(s) >1$ and satisfies the automorphic 
relation 
$$
E^{(n)}(\gamma z,s) = \overline{\kappa(\gamma)}E^{(n)}( z,s).
$$
The series $E^{(n)}( z,s)$ can be expanded in a Fourier series, and the 
constant coefficient is
$$
A_0(s,y) = y^{2s} + \frac{\zeta_K^*(2ns-n)}{\zeta_K^*(2ns-n+1)}y^{2-2s},
$$
where $\zeta_K(2ns-n)$ is the zeta function of the underlying field with
completion $\zeta_K^*$.  This 
has a simple pole when $2ns-n=1$, i.e, at $s = 1/2 +1/(2n)$.   Taking 
the residue at this point, Kubota defined the theta function on the $n$-cover 
of $\GL_2$ by
$$
\theta^{(n)}(z) = \text{\rm Res}_{s= 1/2 +1/(2n)} E^{(n)}( z,s).
$$
Ignoring non-generic primes, the series $E^{(n)}( z,s)$ has a Fourier 
expansion of the form
$$
E^{(n)}( z,s) = A_0(s,y) +y \sum_{m \ne 0}A_m(s)N_{K/\Q}(m)^{s-1/2}K_{2s-1}(4 \pi|m|y)e(mx).
$$
Here $e(x)$ is an additive character with kernel the ring of integers of $K$. 
The coefficients are written as an arithmetic part  multiplied by a $K$-Bessel 
function.  The arithmetic part is 
$$
A_m(s) = \sum_{d \equiv 1 \bmod N} \frac{g_n(m,d)}{{{N_{K/\Q}}}(d)^{2s}}.
 \eqno(1.1)
$$
This is a Dirichlet series built from Gauss sums:
$$
g_n(m,d) = \sum_{r \bmod d}\left(\frac r d\right)_n e\left(\frac{r m}{d} \right).
$$
If we write the Fourier expansion of $\theta^{(n)}(z)$ as
$$
\theta^{(n)}(z) = \tau_n(0)y^{1-1/n}+ y\sum_{m\ne 0}\tau_n(m)K_{1/n}(4 \pi |m|y),
$$
then
$$
\tau_n(m)=N_{K/\Q}(m)^{1/(2n)} \text{\rm Res}_{2s=1+1/n}A_m(s).
$$
The question facing us is the determination of the nature of the 
coefficients $\tau_n(m)$. In the remainder of this section, we normalize
the Fourier expansion to have $\tau_n(1) = 1$.

The Gauss sums factor in the following way: if $d = d_1d_2$, 
with $(d_1,d_2)=1$, then
$$
g_n(m,d_1d_2)= g_n(m,d_1)g_n(m,d_2)\left(\frac {d_1}{ d_2}\right)_n\left(\frac {d_2}{ d_1}\right)_n.
$$
Thus, if $n=2$, the two quadratic symbols cancel at all but finitely many 
places, and the Dirichlet series (1.1) factors into an Euler product which, 
up to a finite number of factors, equals $L_K(2s-1/2,\chi_m)$,  the  
Hecke $L$-series associated to  the quadratic extension $K(\sqrt{m})$.  This 
has a pole at $s=3/4$ when $m$ is a square, explaining why the residue 
of $E^{(2)}( z,s)$, which is the quadratic theta function over the 
field $K$, has a Fourier expansion supported by the square indices. 

For $n \ge 3$ the product $\left(\frac {d_1}{ d_2}\right)_n
\left(\frac {d_2}{ d_1}\right)_n$ is not trivial, and the Dirichlet 
series (1.1) does not factor into an Euler product.   This has so far made it 
impossible to analyze $A_m(s)$ and compute its residue directly.   
Patterson, [\PatCub], was able to use a converse theorem to show that in 
the case $n=3$ the Mellin transform of $\theta^{(n)}(z)$ essentially 
equaled $A_m(s)$,  the first Fourier coefficient of $E^{(3)}( z,s)$.   As 
a consequence, he discovered that the coefficients $\tau_3(m)$ satisfy a 
periodicity relation:
$$
\tau_3(m^3d)= N_{K/\Q}(m)^{1/2}\tau_3(d).
$$
Also, for $d$ cube free, $\tau_3(d) =0$ if $p^2\mid d$ for any prime $p$, and 
for $d$ square free,
$$
\tau_3(d) = \frac{\overline{g^{(3)}(1,d)}}{ N_{K/\Q}(d)^{1/2}}.
$$
Suzuki, [\Suz], attempted to generalize Patterson's method to $n=4$, but 
only succeeded in obtaining partial information about the $\tau_4(m)$.  
Deligne, studying this problem from a representation theoretic point of view, 
was able to explain that  the inaccessibility  of the cases $n \ge 4$ was 
due to a phenomenon of non-uniqueness of Whittaker models.   This approach was 
greatly generalized in a paper of Kazhdan and Patterson [\KP].   In this 
paper, they showed that the periodicity property held in great generality.   
They also showed that the theta functions were eigenfunctions of Hecke 
operators, and that a certain subset of the coefficients were determined by 
these operators.   In the case $n=3$ this subset was everything, but 
for all $n\ge 4$ the coefficients were only partially determined.  

For each prime $p$ there is an associated Hecke operator  $T_{p^n}$.  The eigenvalue of $\theta^{(n)}(z)$ is 
$$
\lambda_{p^n}= N_{K/\Q}(p)^{1/2}+ N_{K/\Q}(p)^{-1/2}.   
$$
To describe the effect of $T_{p^n}$ it will be useful to introduce the 
following notation.  For $0 \le j \le n-1$,
$$
G_j(m,d) =  \frac{\sum_{r \bmod d}\left(\frac r d\right)_n^j e\left(\frac{r m}{d} \right)}{ N_{K/\Q}(d)^{1/2}}.
$$
This is simply the Gauss sum with numerator $m$ and denominator $d$, formed 
with the $j^{th}$ power of the residue symbol, and normalized to have 
absolute value 1 when $d$ is square free.

Applying $T_{p^n}$  to $\theta^{(n)}(z)$ forces the following relation upon 
the coefficients $\tau_n(m)$.  For $(m,p)=1$,
$$
\lambda_{p^n}\tau_n(mp^j)= \tau_n(mp^{j+n})+\tau_n(mp^{j-n}) +  
           N_{K/\Q}(p)^{-1/2}G_{j+1}(m,p)\tau_n(mp^{n-2-j}).
$$
We adopt the convention that $\tau_n(a)$ vanishes unless $a$ is an integer.
The periodicity established in this context in [\Hoff, \KP] is
$$
\tau_n(mp^n)=\tau_n(m) N_{K/\Q}(p)^{1/2}.
$$
For $j = n-1$, the above becomes
$$
( N_{K/\Q}(p)^{1/2}+ N_{K/\Q}(p)^{-1/2})\tau_n(mp^{n-1})= \tau_n(mp^{{n-1}}) N_{K/\Q}(p)^{1/2},
$$
which forces $\tau_n(mp^{n-1})=0$.   For $0 \le j \le n-2$, we obtain
$$
\tau_n(mp^{j})= G_{j+1}(m,p)\tau_n(mp^{n-2-j}).
$$
In the case $n=2$, this means that we have $\tau_2(m)=0$ if $m$ is 
not a square, and $\tau_2(m^2)=N_{K/\Q}(m)^{1/2}$,
a complete description of $\tau_2(m)$.
When $n=3$, we see that $\tau_3(mp^{2})=0$, and 
$$
\tau_3(mp^{}) =  G_{2}(m,p)\tau_3(m).
$$
Reminding our normalisation $\tau_3(1) = 1$, this yields 
$$
\tau_3(cd^3) = N_{K/\Q}(d)^{1/2}\overline{G_1(1,c)},
$$
for $c$ square free, and $\tau_3(m) =0$ otherwise.  This is a complete 
description of $\theta^{(3)}(z)$, which agrees with that found by Patterson.

When $n=4$, the first example of undetermined coefficients occurs.   We see 
that for $(m,p)=1$, $\tau_4(mp^3)=0$.  Also 
$$
\tau_4(mp^{2}) =  G_{3}(m,p)\tau_4(m)
$$
and
$$
\tau_4(mp) =  G_{2}(m,p)\tau_4(mp).
 \eqno(1.2)
$$
Taking $m=1$, we see that although $\tau_4(p^{2})$ is 
determined, $\tau_4(p^{})$ is not.
Interestingly, as the quadratic Gauss sum is trivial in this context, we have
$$
G_2(m,p) = \left(\frac m p\right)_4^2 =   \left(\frac m p\right)_2.
$$
It follows then, from the above, that if $ \left(\frac m p\right)_2 = -1$, 
then $\tau_4(mp)=0$.   More generally, if $m$ possesses any 
factorization $m = m_1m_2$, with $ \left(\frac {m_1}{ m_2}\right)_2 = -1$,
then $\tau_4(m)=0$. 

When $n=5$, one finds that $\tau_5(p^4)  =0$, $\tau_5(p^3)  
               =\overline{G_1(1,p)}$, and that
$$
\tau_5(p) = G_2(1,p)\tau_5(p^2).
$$
This finally leads us to the subject of this paper.   When $n=6$, the Hecke 
relations imply that $\tau_6(p^5)  =0$, $\tau_6(p^4)  =\overline{G_1(1,p)}$, 
that
$$
\tau_6(p) = G_2(1,p)\tau_6(p^3), 
$$
and that $\tau_6(p^2)$ is related to itself via
$$
\tau_6(mp^{2})= G_{3}(m,p)\tau_6(mp^{2}).
 \eqno(1.3)
$$
Interestingly, the Gauss sum appearing  in (1.3) is quadratic, as in (1.2), 
suggesting a possible parallel phenomenon occurring in the cases $n=4$ 
and $n=6$. We will see in Section~5 that the relation 
$\tau_6(p^4) = \overline{G_1(1,p)}$ almost holds in a more precise setup.

What rule or pattern, if any, governs the undetermined coefficients?  One 
striking observation and conjecture was made by Patterson in the case $n=4$.   
Recall that the first Fourier coefficient of 
$E^{(4)}( z,s)$ was
$$
A_1(s)= \sum\frac{G_1(1,m)}{N_{K/\Q}(m)^{2s-1/2}}.
$$
As $A_1(s)$ is a Fourier coefficient of $E^{(4)}( z,s)$, which possesses a 
functional equation as $s \rightarrow 1-s$, $A_1(s)$ inherits the same 
functional equation.  Change the variable, rename this series as
$$
\psi(w) = \sum\frac{G_1(1,m)}{N_{K/\Q}(m)^{w}},
$$
and consider the Dirichlet series $D_1(w) = \zeta_K(4w-1)\psi(w)$.   This 
has a functional equation as $w \rightarrow 1-w$, and a simple pole 
at $w = 3/4$.  On the other hand, the Dirichlet series 
$$
D_2(w) = \zeta_K(4w-1)\ \sum\frac{\tau_4(m)^2}{N_{K/\Q}(m)^w}
$$
is the Rankin-Selberg convolution of $\theta^{(4)}(z)$ with itself and can be 
easily seen to have a functional equation as $w \rightarrow 1-w$, and a 
double pole at $w = 3/4$.   Patterson observed that the gamma factors 
occurring in the functional equations of $D_1(w)^2$ and $D_2(w)$ were 
identical, and conjectured that
$$
\overline{D_1(w)}^2=D_2(w).
$$
This conjectured equality can be seen to be consistent with all the 
information provided by periodicity and the Hecke relations.   Dividing by an 
extra $\zeta_K(4w-1)$, the conjecture states that
$$
\sum\frac{\tau_4(m)^2}{N_{K/\Q}(m)^w}= \zeta_K(4w-1) \left(  \sum\frac{\overline{G_1(1,m)}}{N_{K/\Q}(m)^{w}} \right)^2.
$$
In other words, the conjecture predicts the values of $\tau_4(m)$ up to sign.  
Checking the coefficients of $m=p^2$, we see that on the left hand side we 
have $\tau_4(p^2)^2= \overline{G_1(1,p)}^2$, while on the right hand side, 
as $G_1(1,p^2) =0$, the only contribution comes from the square of the $m=p$ 
term, namely  $\overline{G_1(1,p)}^2$.  Checking further, for $m$ square 
free, on the right hand side we have
$$
\sum_{m = m_1m_2}\overline{G_1(1,m_1)}\overline{G_1(1,m_2)}= \overline{G_1(1,m)}\sum_{m = m_1m_2}\left(\frac{m_1}{m_2}\right)_2,
$$
which does indeed vanish  if $m$ possesses any factorization $m = m_1m_2$, 
with $ \left(\frac {m_1}{ m_2}\right)_2$ $= -1$.   Most interestingly, looking 
at the prime indices, the conjecture predicts that 
$$
\tau_4(p)^2 = 2\overline{ G_1(1,p)}.
$$
In [\Var] a conjecture was made about the $n=6$ case that was weaker than 
the $n=4$ conjecture, in that it did not quite pin down all of the 
coefficients.   This conjecture was that
$$
\sum \frac{\tau_6(m^2)}{Nm^u} =\sum  \frac{\overline{\tau_3(m)}}{Nm^u}
\cdot  \sum \frac{G_1^{(3)}(1,d)}{Nd^u},
$$
where the superscript $(3)$ indicates that we considering the {\it cubic\/}
Gauss sum.
The left hand side is the convolution of the theta function on the $6$-cover 
of $\GL_2(\C)$ with the theta function on the $2$-cover of $\GL_2(\C)$.    This 
has the effect of picking off Fourier coefficients with square indices.  The 
right hand side is the product of the Mellin transform of the theta function 
on the 3-cover of $\GL_2(\C)$, with the first coefficient of the cubic 
Eisenstein series.  The two, however, are equal in this cubic case, up to a 
zeta function factor.  Writing $m = m_1m_2^2m_3^3$, with $m_1,m_2$  square 
free and relatively prime, $m_3$ unrestricted we see by the periodicity 
properties of $\tau_6$ and the known valuation of $\tau_3$ that 
after canceling a zeta factor on both sides this relation translates to 
$$
\sum \frac{\tau_6(m_1^2m_2^4)}{Nm_1^uNm_2^{2u}} =
\left(\sum \frac{G_1^{(3)}(1,d)}{Nd^u}\right)^2,
$$
another curious identity involving the square of a series without an Euler 
product.  Note that the Gauss sums $G_1^{(3)}(1,d)$ on the right hand side 
vanish unless $d$ is square free.

Equating corresponding coefficients we have the  following predicted behavior 
for the coefficients $\tau_6(m_1^2m_2^4)$:
$$
\tau_6(m_1^2 m_2^4) = G_1^{(3)}(1,m_2)^2G_1^{(3)}(1,m_1)\left(\frac{m_2}{m_1}\right)_3^2
\sum_{m_1 = d_1d_2}\left(\frac{d_2}{d_1}\right)_3.
$$
In particular, when $m_1=p$ and $m_2=1$, this reduces to the 
relation $\tau_6(p^2) = 2 G_1^{(3)}(1,p)$.  This the fundamental 
relation which is being tested in this paper.

We will see in Section~5 that computational evidence overwhelmingly supports
the conjecture for $|p| \equiv 7 \bmod 12$. Indeed, our computations suggest
that, apart from a $12$-th root of unity, 
$$
\tau_6(p^2) = 2 {G_1^{(3)}(1,p) \over \sqrt{3}}.
$$
We recall the conjecture from [\Var] was made disregarding
the prime 3, so it should come as no surprise that an additional power of~$3$
occurs in the actual coefficients. We remark that care should be made in
comparing the current article and [\Var], since the definition of the sixth
order symbol in the two articles are {\it conjugates\/} of each other. 

We will give a conjecture for $\tau(p^2)$ for the other congruence classes
of $|p|$ in Section~5. We give a proof of certain special cases as well. 
Finally, we examine the square~$\tau(p)^2$ and give a conjectured value for
this coefficient.

\noindent
\head 2. Theta series
\endhead
\noindent
Throughout this section, we fix an integer~$n>2$. We let $K = \Q(\zeta_n)$
be the cyclotomic field obtained by adjoining a primitive $n$-th root of
unity~$\zeta_n$. Later on, we will focus on~$n=6$, and to make the
exposition easier, we restrict ourselves to the case that $K$ has
class number one in this article. We define the set
$$
S = \{ v_\pi \hbox{\ with\ } \pi \mid n\infty\}
$$
to be the places dividing~$n$ together with the infinite places. We note that
since~$K$ is totally imaginary, all infinite places are complex. The
set of all finite places dividing~$n$ is denoted by~$S_f$. We let
$$
K_S = \prod_{v\in S} K_v
$$
be the product of the completions at all the places in~$S$. We embed $K$ into the product $K_S$
along the diagonal. Our first goal in this section is
to define a {\it Gauss sum\/} on the ring 
of $S$-integers~$\O_S = \O[\pi^{-1} \mid v_{\pi} \in S_f]$.

\subhead 2.1 Gauss sums
\endsubhead

For $v\in S$, the localization $K_v$ admits a generalized Hilbert symbol. We
recall its construction here. We let $L$ be a local field of characteristic
zero with $\zeta_n \in L$, and we let $M = L(\root n \of{L^*})$. By local Artin
reciprocity, we have 
$$
L^* / N_{M/L}(M^*) \cong \Gal(M/L)
$$
via the Artin map. The equality $L^{*n} = N_{M/L}(M^*)$ and Kummer theory 
give a map
$$
L^*/L^{*n} \cong \Hom(\Gal(M/K),\mu_n).
$$
Combining both displayed equations gives the {\it Hilbert symbol\/}
$$
(x,y): L/L^{*n} \times L/L^{*n} \rightarrow \mu_n
$$
as $(x,y) = \chi_y((x,M/L))$. Here, $\chi_y$ is the Kummer character of~$y$,
and $(\cdot,M/L)$ is the Artin symbol. We combine the local Hilbert symbols
to a symbol on $K_S$ via
$$
(x,y)_S = \prod_{v \in S} (x,y)_v.
$$

For coprime $a,b \in \O_S$, we let 
$$
{a \legendre b}_S = \prod_{v \not \in S, v \mid b} (a,b)_v
$$
be the generalized Legende symbol. The Hilbert symbol and the Legendre
symbol satisfy a {\it reciprocity law\/}
$$
{a \legendre b}_S = (a,b)_S {b \legendre a}_S
$$
that will be useful for explicit computations, see Section~3. We note 
that the Hilbert symbol is local and defined for all $a,b\in K_S$, whereas the 
Legendre symbol is global and more restricted.

Having defined a multiplicative character on $K_S$, we now proceed with
defining an additive character~$e$. As before, we will do so by defining a local character $e_v$ for each~$v\in S$. The 
desired character~$e$ is then simply the product of the local characters. 
First assume that $L$ is $\gp$-adic, i.e., a finite extension of the $p$-adic 
field~$\Q_p$.  We define $e_v$ as the composition of the maps
$$
L \mapright{\Tr} \Q_p \mapright{} \Q_p/\Z_p \mapright{\lambda} \Q/\Z \mapright{\exp(-2\pi i \cdot)} \C^*.
$$
Here, the map $\lambda: \Q_p \mapsto \Q/\Z$ satisfies
$$
\lambda\left(\sum_{j=-N}^{\infty} x_j p^j\right) = \sum_{j=-N}^{-1} x_j p^j,
$$
i.e., it is the `tail' of the $p$-adic extension of~$x$, so that 
$\Ker(\lambda) = \Z_p$.  For $L=\C$, we 
put $e_\infty(x) = \exp(-2\pi i \Tr(x))$.

For an embedding $\varepsilon: \mu_n(K) \rightarrow \C$, we can now define
the {\it Gauss sum\/}
$$
g_n(r,\varepsilon,c) = \sum_{x \bmod c} \varepsilon\left( {x \legendre c}_S\right) e(rx/c)
$$
for $r,c \in \O_S$. We will present an algorithm to compute~$g_n(r,\varepsilon,c)$ in
Section~3.

\subhead 2.2. Dirichlet series
\endsubhead
\ \par\noindent
The Dirichlet series we will be working with are 
indexed by the group~$K_S^*/(K_S^{*n}\O_{S}^*)$.
We first explain the structure of this group.

\proclaim{2.1 Lemma}
We have $\O_{S}^* \cap K_S^{*n} = \O_{S}^{*n}$.
\endproclaim
\noindent{\bf Proof.} Let $x \in \O_{S}^* \cap K_S^{*n}$. Since $x$ is
locally an $n$-th power, we have $(x,c)_S = 1$ for all $c \in \O_{S}$
by the properties of the Hilbert symbol.
This means that ${x \legendre c}_S = 1$ holds, which implies that $x$ arises
from a global $n$-th power. Hence, $x \in \O_{S}^{*n}$.\hfill\endproof

\proclaim{2.2 Lemma}
The following equality holds:
$$
[K_S^* : K_S^{*n}\O_{S}^*] = n^{\#S}
$$
\endproclaim
\noindent
{\bf Proof.} This is proven in~[\PatIJM, Sec.\ 3]. We give a slightly
modified proof here for convenience. 
By standard group theory, we have
$$
[K_S^* : K_S^{*n}\O_{K,S}^*] = [K_S^* : K_S^{*n}]/[\O_{S}^*K_S^{*n} : K_S^{*n}]
                             = [K_S^* : K_S^{*n}]/[\O_S^* : \O_S^* \cap K_S^{*n}],
$$
and by Lemma~2.1 we have $[\O_S^* : \O_S^* \cap K_S^{*n}] = [\O_S^* :
\O_S^{*n}]$. Using Dirichlet's unit theorem, we compute this last index to
be~$n^{\#S}$. It remains to compute $[K_S^* : K_S^{*n}]$. For a finite
place $v_p \in S$, we have
$$
[ K_{v_p}^* : K_{v_p}^{*n} ] = {n^2 \over p^{-v_p(n)}}
$$
by [\Neu, Cor.~II.5.8].  We let $S_\infty$ be the set of infinite 
places in~$S$.  By the product formula~[\Neu, Prop.~III.1.3], we have
$$
\prod_{v_p \in S_f} p^{-v_p(n)} = \left(\prod_{v_p \in S_\infty} N_{\C/\R}(n)
\right)^{-1} = n^{-2\#S_\infty},
$$
and we conclude that we have $[K_S^* : K_S^{*n}] = n^{2\#S}$. The lemma 
follows.  \hfill\endproof\par
\ \par\noindent
We pick a coset $\eta$ of $K_S^*/K_S^{*n}$. 
For $r \in K_S$, and $s \in \C$ with $\text{\rm Re}(s) > 1/2$, we define the 
Dirichlet series
$$
\psi(r,s,\eta) = \sum_{c \in (\eta K_S^{*n} \cap \O_S)/\O_S^{*n}}
g_n(r,\varepsilon,c) |c|_S^{-s-1} L_S(|\cdot|_S^{ns+1}),
$$
where 
$$
L_S(|\cdot|_S^s) = \prod_{v\in S_f}(1-|\pi_v|_v^s)\zeta_K(s).
$$
In the last expression, $\pi_v$ is a uniformizer for~$K_v$. The norm $|c|_S$
appearing in $\psi$ is the `$S$-norm' $|c|_S = \prod_{v\in S} |c_v|_v$. The
$S$-norm coincides with the regular norm $N_{K/\Q}$ for $\gcd(c,n)=1$.

\proclaim{2.3. Remark}
For $u \in \O_S^*$, we have
$$
\psi(r,s,u\eta) = \varepsilon(u,\eta)_S \psi(r,s,\eta).
$$
\endproclaim

From Remark~2.3, we see that for understanding $\psi(r,s,\eta)$ it 
suffices to pick a coset $\eta$ for $K_S^*/(K_S^{*n}\O_S^*)$. The following
theorem describes the analytic continuation of~$\psi$.

\proclaim{2.4. Theorem}
The function $\psi$ defined above converges absolutely 
for~${\text Re}(s) > 1/2$. It admits a memorphic extension to~$\C$. This
extension has for~${\text Re}(s) \geq 0$ at most a single pole at~$s=1/n$. 
\endproclaim
\noindent{\bf Proof.} See [\EP, Sec.~5]. \hfill\endproof\par
\ \par\noindent
The residue at $s=1/n$ of $\psi(r,s,\eta)$ is related to the Fourier 
coefficient~$\tau(r,V)$ from the introduction in the following way. We
pick a full set of cosets $V$ for $K_S^*/(K_S^{*n}\O_S^*)$; this set has 
cardinality~$n^{\#S}$ by Lemma~2.2. After possibly multiplying by some
element of $\O_S^*$, we assume that all $\eta = (\eta_1,\ldots,\eta_{\#S})\in V$
are integral at each component.

We now look at the series
$$
A(r,s,V) = \sum_{\eta\in V} \psi(r,s,\eta) = \sum_{c \in \O/\O^{*n} \atop\scriptscriptstyle c=(c_1,
\ldots,c_{\#S}) \in V} g_n(r,\varepsilon,c) |c|^{-s-1} L_S(|\cdot|^{ns+1}).
$$
The series~$A$ has a pole at $s=1/n$. By comparing the formula above 
with~(1.1), we see that this is in fact the same place as the pole for
$A_m(s)$ from the introduction. 

The last sum in the equation above is a sum over {\it ideals\/}~$I \subset
\O_K$ coprime to~$n$ with the convention that we pick a generator~$c$ 
of~$I$ with $c \in \eta k_S^{*n}$ for some $\eta\in V$. The quantity
$$
\tau(r,V) = N_{K/\Q}(r)^{1/2n}{\text\rm Res}_{s=1/n} A(r,s,V)
$$
is the main object of study in this paper.
Different choices for~$V$ yield different Fourier 
coefficients~$\tau(r,V)$. The introduction takes $d \equiv 1 \bmod N$, but
there are other choices one can make.  
As we will see in Section~5, selecting a convenient~$V$ is
part of our conjecture.

\noindent
\head 3. Computing Fourier coefficients
\endhead
\noindent
The functions~$\psi(r,s,\eta)$ satisfy a functional equation 
in~$s \rightarrow -s$. To state the equation, we modify $\psi$ slightly
and define
$$
\Psi(r,s,\eta) = y_1^{-s} \left( {\Gamma(ns)\over\Gamma(s)} \right)^{N/2} 
\psi(r,s,\eta),
$$
where $N = [K:\Q]$ and 
$$
y_1 = (2\pi)^{N(n-1)\over 2}|r|_S^{-1/2}\prod_{v \in S_f} {|\pi_v|}_v^{nd_v/2},
$$
with $d_v$ the local different as in Section~3.1.
\proclaim {3.1. Functional Equation} [\EP, Sec.~5] The function~$\Psi$ satisfies a
functional equation
$$
\Psi(r,s,\eta_i) = n^{N/2} G_f(s) \prod_{v\in S_f } |\pi_v|_v^{-d_v/2} \sum_{j=1}^{n^{\#S}} T_{ij}(r,s) \Psi(r,-s,\eta_j), \eqno(3.1)
$$
where $\{ \eta_i \}_i$ is a full set of representatives for $K_S^*/(K_S^{*n}\O_S^*)$.
\endproclaim
Before we explain the notation in equation (3.1) above, we note that it
translates a {\it vector\/} $\Psi(r,s,\eta_i)_i$ to $\Psi(r,-s,\eta_i)_i$.
This means that although we are interested in computing a sum over 
all $\eta_i$, we do have to work with the individual~$\Psi(r,s,\eta_i)$'s.

In the functional equation, we have
$$
G_f(s) = \prod_{v \in S_f} {1 \over 1-|\pi_v|_v^{1-ns}}.
$$
The integer $d_v$ for a place $v$ over~$p$ is related to 
the {\it different\/} of the extension $\Q_p(\zeta_n)/\Q_p$ in the 
following way. The maximal order of~$\Q_p(\zeta_n)$ equals $\Z_p[\zeta_n]$
and we let $f\in \Z_p[x]$ be the minimal polynomial of~$\zeta_n$. The
different of $\Q_p(\zeta_n)/\Q_p$ equals $(f'(\zeta_n))$ and we have
$$
(f'(\zeta_n)) = (\pi_v)^{d_v}.
$$
In particular, $d_v = 0$ if $\Q_p(\zeta_n)$ is unramified.

Finally, the $T_{ij}(r,s)$ occuring in (3.1) are 
coefficients of an $n^{\# S}\times n^{\# S}$-matrix $T(r,s)$. The matrix~$T$ 
is defined over~$\overline\Q (|\pi_v|_v^s \mid v \in S_f)$. 
We will give a method to compute~$T_{ij}(r,s)$ below. 

Knowing the functional equation that $\Psi(r,s,\eta)$ satisfies, the idea
is to compute its residue at $s = 1/n$ using contour integration. Indeed, 
since the function $\Psi(r,s,\eta)$ decays exponentially fast for Im$(s) \rightarrow
\infty$ by [\EDiss, Thm.~2.4], we have
$$
\text{\rm Res}_{s=1/n} \Psi(r,s,\eta) = {1 \over 2\pi i} \int_{(\sigma)} \Psi(r,s,\eta)\text{d}s - {1 \over 2 \pi i}\int_{(-\sigma)} \Psi(r,s,\eta)\text{d}s, \eqno(3.2)
$$
where $(\sigma)$ denotes the vertical line Re$(s) = \sigma > 1/2$. It will
be convenient for our computations to modify (3.2) slightly. If $f$ is 
a holomorphic function such that $f(s) (\Gamma(ns)/\Gamma(s))^{N/2}$ decays 
exponentially for Im$(s) \rightarrow\infty$, then we also have
$$
f(1/n)x^{-1/n}\text{\rm Res}_{s=1/n} \Psi(r,s,\eta) = {1 \over 2\pi i}\int_{(\sigma)-(-\sigma)} \Psi(r,s,\eta) f(s) x^{-s} \text{d}s \eqno(3.3),
$$
with $(\sigma)-(-\sigma)$ the union of the two vertical lines. Not only
does (3.3) allow for greater flexibility in computing the two integrals, but
it also servers as a check on our computations by letting the parameter~$x>0$
vary.

\subhead 3.1. Integral with $\sigma>0$
\endsubhead
\par \noindent
In this subsection we explain how to approximate
$$
{1\over 2\pi i} \sum_{\eta\in V} \int_{(\sigma)}\Psi(r,s,\eta) f(s) x^{-s} \text{d}s \eqno(3.5)
$$
for $\sigma>1/2$. Using the formula for $A(r,s,V)$, we are interested in
computing the sum
$$
{1\over 2\pi i}\sum_{c\in \O/\O^{*n} \atop \scriptscriptstyle c \in V} \int_{(\sigma)}{g_n(r,\varepsilon,c)\over
|c|} (y_1|c|x)^{-s} L_S(|\cdot|^{ns+1}) \left({\Gamma(ns)\over\Gamma(s)}\right)^{N/2} f(s) \text{d}s.
$$
In order to compute this sum, we first write $L_S(|\cdot|^{ns+1}) 
= \sum_{m=1}^{\infty} a(m) m^{-s}$. Since $L_S$ is basically the Dedekind
zeta-function of~$K$, the coefficients~$a(m)$ are easily determined. Indeed,
we have
$$
a(m) = \cases I(k)/k, & \text{if\ } m=k^n \text{\ and if\ } m_v = 1 
\text{\ for\ } v\in S_f, \\
0, & \text{otherwise,}\\ 
\endcases\eqno(3.6)
$$
with $I(k)$ the number of $\O$-ideals of norm~$k$.
We note that the coefficients~$a(m)$ decay quite rapidly since they
are only supported on~$n$-th powers. If we now put
$$
F_1(x) = {1 \over 2\pi i} \int_{(\sigma)} \left( {\Gamma(ns)\over\Gamma(s)}
\right)^{N/2} f(s) x^{-s} \text{d}s
$$
then the sum in (3.5) is equal to 
$$
\sum_{c\in\O/O^{*n} \atop \scriptscriptstyle c \in V} {g_n(r,\varepsilon,c) \over
|c|} \sum_{m=1}^{\infty} a(m) F_1(x y_1 m |c|). \eqno(3.7)
$$
A key idea to approximating the sum above is to loop over all {\it ideals\/}
of $\O$ that are coprime to~$n$. For every ideal $I$, we only
consider the generators $c$ that lie in~$V$ when viewed as elements of $K_S$. 
If we can efficiently compute
the function~$F_1$, then for every such~$c$, we 
approximate the sum~$\sum_m a(m) F_1(xy_1 m|c|)$ to high precision and get
the contribution coming from~$c$ to the sum~(3.7). 

The inclusion of the function~$f(s)$ in the integral gives us many choices
for the function~$F_1(x)$. There is a trade off in picking $f$ so that $F_1$
converges fast, and is easy to evaluate. We refer to Section~4 for an 
example.

\subhead 3.1. Integral with $\sigma<0$
\endsubhead
\par \noindent
The main idea behind evaluating (3.5) for $\sigma<0$ is the same as 
for $\sigma>0$. The only technical difficulty is that since $\Psi$ does not 
admit a Dirichlet expansion for $\sigma < 1/2$, we will map $s$ to $-s$ and
use the functional equation. After replacing $s$ by $-s$ in expression~(3.5),
we see that we have to evaluate
$$
{1 \over 2\pi i} \sum_{\eta_i\in V} \int_{(-\sigma)} G_f(-s) f(-s) \left({1\over x}\right)^{-s} \left( \sum_{j=1}^{n^{\#S}} T_{ij}(r,-s) \Psi(r,s,\eta_j)\right) \text{d}s. \eqno(3.8)
$$
We note that $G_f(-s) L_S(|\cdot|^{ns+1}) = \zeta_K(ns+1) = \sum_m b(m) m^{-s}$
holds. The coeffcients $b(m)$ satisfy $b(m) = I(k)/k$ if $m = k^n$, and
$b(m) = 0$ otherwise. Just like in the previous subsection, these coefficients
decay rapidly. Similar to the case $\sigma>0$, we put
$$
F_2(x) = {1 \over 2\pi i} \int_{(\sigma)} \left( {\Gamma(ns)\over\Gamma(s)}
\right)^{N/2} f(-s) x^{-s} \text{d}s.
$$
The only difference with the previous subsection is that we have to deal with
the sum $\sum_i T_{ij}(r,-s)$ of the coefficients in the $j$-th column of 
the transition matrix. 
Because each coefficient $T_{ij}$ lives in $\overline\Q(|\pi_v|_v^s \mid
v \in S_f)$, we can write
$$
\sum_{i=1}^{n^{\# S}} T_{ij}(r,-s) = \sum_{w \in \Z^{\# S_f}} c_{jw} |\pi|^{ws}
\qquad \text{with} \qquad |\pi|^w = \prod_{v \in S_f} |\pi_v|_v
$$
for certain algebraic numbers~$c_{jw}$. Putting everything together, we
can write
$$
{1\over 2\pi i} \sum_{\eta\in V} \int_{(\sigma)}\Psi(r,s,\eta) f(s) x^{-s} \text{d}s =
$$
$$
= y_2 \sum_{c \in \O/\O^{*n} \atop \scriptscriptstyle c \in V} {g_n(r,\varepsilon,c) \over |c|} \sum_{w \in \Z^{\#S_f}} c_{jw} \sum_{m=1}^{\infty} 
b(m) F_2(x^{-1} y_1 |c| |\pi_v|_v^w m), \eqno(3.9)
$$
with $y_2 = n^{N/2} \prod_{v} |\pi_v|_v^{d_v/2}$.

Just as for the integral with $\sigma>0$, we will loop over all {\it ideals\/}
of $\O$ that are coprime to~$n$ and for every ideal $I$, we only
consider the generators $c$ that lie in~$V$ when viewed as elements of~$K_S$. 

We observe that for each $c\in\O$, we need to evaluate the function~$F_2$
many more times than~$F_1$. Indeed, whereas we only need to evaluate $F_1$
once for every~$m \geq 1$, we need to evaluate $F_2$ for every $w$ such that
$c_{jw}$ is nonzero. This means that we need to pick our function~$f$ so that
$F_2$ is particulary easy to evaluate. 

\subhead 3.1. The transition matrix
\endsubhead
\par \noindent
We recall that for every choice $V$ of cosets for $K_S^*/(K_S^{*n}\O_S^*)$,
there exists a matrix
$$
T(r,s) \in \overline\Q(|\pi_v|_v^s \mid v \in S_f)
$$
such that $\Psi(r,s,\eta_i)_{\eta_i \in V}$ satisfies functional
equation~(3.1). In this subsection we explain how to compute the coefficients
of this matrix. We will only give the results needed for actual computations,
and refer to [\EP, Sec.~5] for the underlying theory.

We fix a choice of coset representatives~$V$. One can show that $T_{ij}$
satisfies
$$
T_{ij}(s,r) = \left({r \over \eta_i\eta_j}\right)_{S_f}^s {(\eta_i,-\eta_j)_S\over
n^{2\#S}} \prod_{v \in S_f} (1-|\pi_v|_v^{ns}) \cdot
$$
$$
\cdot \sum_{h \in \O_S/\O_S^{*n}} \left(h, -{\eta_j\over\eta_i}\right)_S
\prod_{v \in S_f} \sum_{y \in K_v^*/K_v^{*n}} \left(-{r_v \over h
\eta_i\eta_j},y\right)_S \Gamma_{v}(|\cdot|_v^s\varepsilon(y,\cdot)).
\eqno(3.10)
$$
\hfuzz=1.5pt
As indicated in the previous subsection, we will view $|\pi_v|_v$ as 
indeterminates and compute the coefficients~$T_{ij}$ as elements of an
$\#S_f$-dimensional function field over~$\overline\Q$. To make the
computations as fast as possible, we should be careful to select the
number field~$L$ over which the coefficients of~$T_{ij}$ are defined. We
will see in the next section that we can take the Euclidean 
field $L = \Q(\zeta_{36})$ for $n=6$ for instance. 

Before we detail the computation of the local Gamma function~$\Gamma_{v}$,
we explain an idea from [\Well] to speed up the computations of all 
the~$T_{ij}$. By examining~(3.10) closely, we see that except for a 
factor~$(\eta_i,-\eta_j)_S$, it only depends on $\eta_i^2$ and
$\eta_i\eta_j$. Hence, for even~$n$, we can save time by only computing
the coefficients in~(3.10) for $\eta_i^2$ and $\eta_j$. 

The factor~$\Gamma_v$ occuring in~(3.10) was introduced by Tate in his
thesis [\Tate], where it is called~$\rho(c)$. We recall some of the basic
theory here. We let $\chi: K_v^* \rightarrow
\C^*$ be a quasi-character, i.e., a continuous multiplicative (but not
necessarily of absolute value 1) map. We have $K_v^* \cong (\pi_v) \times U$,
and $\chi$ equals $\widetilde\chi |\cdot|^s$ by [\Tate, Thm.\ 2.3.1], 
where $\widetilde\chi$ is a {\it character\/} on the unit group~$U$ of
the maximal order of~$K_v$. We furthermore have
$U = \mu_{q-1} \times U^{(1)}$ with $q$ the characteristic of the residue
field, and $U^{(1)} = 1+\gp$. If $\widetilde\chi$ is trivial on~$U$, we
call $\widetilde\chi$ {\it unramified\/}. Otherwise, since the
subgroups $1+\gp^n$ form a filtration of~$U$, there exists a
minimal $f>0$ with with $\widetilde\chi(1+\gp^f) = 1$. In this ramified
case, we call the integer $f$ the {\it conductor\/} of~$\widetilde\chi$.

The character we are interested in for~$(3.10)$ is $\widetilde\chi = \varepsilon(y,\cdot)$. We compute its conductor in the following
way. We have an 
isomorphism $U^{(n)}/U^{(n+1)} \cong \O_v/\pi_v$ for all~$n$, and we 
first compute
a set of representatives~$R_1$ for $\O_v/\pi_v$. We now {\it guess\/} that
$\varepsilon(y,\cdot)$ is trivial on~$\gp^g$ for some~$g$, like $g=n$. We
check that our guess is correct by computing
$$
\varepsilon(y,1+r \pi^g) \qquad \hbox{for all} \qquad r \in R_1.
$$
If the computation above does not yield~1 for all $r$, then we replace
$g$ by $2g$ and repeat the check until we do get 1 for all~$r \in R_1$. If
we do get $1$ for all~$r$, we replace~$g$ by $g-1$ and repeat the check. We
continue doing the latter until we do {\it not\/} get~1. The last~$g>1$ for
which we get 1 for all~$r$ is the conductor of~$\varepsilon(y,\cdot)$. If
we have $\varepsilon(y,1+r\pi) = 1$ for all~$r$, then we do a last check
to see if $\varepsilon(y,x) = 1$ for all $x \in \mu_{q-1}$. If this is
the case, then $\varepsilon(y,\cdot)$ is unramified, otherwise it has
conductor~1.

In the unramified case, we have
$$
\Gamma_v(\chi) = \chi(\pi_v)^{-d} {1-|\pi|_v \chi(\pi_v)^{-1} \over 1-\chi(\pi) } |\pi_v|_v^{d_v/2} \in \overline\Q(|\pi_v|_v^s)
$$
for $\chi = |\cdot|_v^s \varepsilon(y,\cdot)_S$. Here, $d_v$ is as before 
the local different of $\Q_p(\zeta_n)/\Q_p$, and we pick the totally
positive square root of $|\pi_v|_v^{d_v}$. We note that we view $|\pi_v|_v^s$
as an indeterminate, and view the image of $\chi$ in $\overline\Q$. 

In the ramified case, we compute a set of representatives $R_2$ for 
$$
U/(1+\gp^f) \cong (U/\gp^f)^*,
$$
and compute the sum
$$
W(\chi) = |\pi_v|_v^{f/2} \sum_{x \in R_2} \chi(x) e_v\left( {x\over \pi_v^{d+f}}\right)
$$
for $\chi = |\cdot|_v^s \varepsilon(y,\cdot)$. As before, we take the
totally positive square root of~$|\pi_v|_v^f$. The map $e_v$ is the same
character as in Section~2.1, except that we its image in~$\overline\Q$. A good check for the
computations is that $W(\chi)$ is a root of unity. We now have
$$
\Gamma_v(\chi) = \chi(\pi_v)^{-d_v-f} |\pi_v|_v^{{d_v+f}\over 2} W(\chi),
$$
where we again view $|\pi_v|_v$ as an indeterminate.\par
\ \par
\subhead 3.4. Gauss sums
\endsubhead
\par\noindent
The last subsection deals with the computation of the Gauss sum
$$
g_n(r,\varepsilon,c) = \sum_{x \bmod c} \varepsilon\left( {x \legendre c}_S\right) e(rx/c)
$$
that occurs in the sums~(3.7) and~(3.9). The following lemma
reduces the computation to the case $g_n(1,\varepsilon,\pi)$ for a
prime element~$\pi\in\O_S$.

\proclaim{3.2. Lemma} Let $r_,c \in \O_S$ be non-zero with $\gcd(r,c) = 1$. 
Then we have 
$$
g_n(r,\varepsilon,c) = \varepsilon\left( {r \legendre c}_S^{-1}\right)
g_n(1,\varepsilon,c).
$$
Let $c_1,c_2\in  \O_S$ be non-zero with $\gcd(c_1,c_2) = 1$. Then we have
$$
g_n(r,\varepsilon,c_1 c_2) = \varepsilon\left( {c_1\legendre c_2}_S \right)
                             \varepsilon\left( {c_2\legendre c_1}_S \right)
                             g_n(r,\varepsilon,c_1)
                             g_n(r,\varepsilon,c_2).
$$
Let $\pi\in\O_S$ be prime. Then we have
$$
g_n(r,\varepsilon,\pi^k) = 0 \qquad \hbox{for} \qquad \gcd(r,\pi) = 1 \qquad
\hbox{and} \qquad k\geq 2.
$$
Finally, we have
$$
g_n(\pi^k,\varepsilon,\pi^l) = \cases 0 & \hbox{for\ }k \not = l-1\\
N_{K/\Q}(\pi)^k g_n(1,\varepsilon^l,\pi) & \hbox{otherwise,}
\endcases
$$
in the relevant case $\varepsilon^l \not = 1$.
\endproclaim
\noindent{\bf Proof.} See [\PatDis].\  \hfill\endproof\par
\ \par\noindent
We see from this lemma that for computing the residue of~$A(r,s,V)$ we can,
except for the case $\gcd(r,c) \not = 1$, restrict our attention to~$c$
being {\it squarefree\/}. Furthermore, we note that it suffices to make 
a list of Gauss sums $g_n(1,\varepsilon,\pi)$ for prime elements $\pi$ of
norm up to some bound. In the remainder of this section we detail a method
to compute~$g_n(1,\varepsilon,\pi)$ in many cases.

The Gauss sum~$g_n(1,\varepsilon,\pi)$ is closely related to the `ordinary'
Gauss sum
$$
g(\chi) = \sum_{x=1}^{p-1} \chi(x) \exp(2\pi i x/p) \eqno(3.11)
$$
in case~$\pi$ has prime norm~$p$. Here, $\chi$ is the 
character $\varepsilon\left( {x\legendre \pi}_S \right)$ and the $\pi$ in
the exponential is of course the complex number~$\pi\approx 3.14$. The
exact relationship between $g_n(1,\varepsilon,\pi)$ and $g(\chi)$ depends
on $n$, we refer to Section~4 for the case~$n=6$.

The naive way of computing~(3.11) by evaluating the sum directly takes
$\widetilde O(p)$ operations, and this run time can be a bottleneck for the
computations. It is well known that we can do much better, at least
heuristically. We define the {\it root number\/} $W_D(\chi)$ as
$$
W_D(\chi) = \cases {g(\chi) \over m^{1/2}} & \hbox{if\ } \chi(-1) = 1,\\
                   {g(\chi) \over m^{1/2}i}, & \hbox{if\ } \chi(-1) = -1,
\endcases
$$
where we remark that the subscript $D$ serves to distinguish $W_D(\chi)$ from
the root number in subsection~(3.3). We define $e \in \{0,1\}$ so that
$\chi(-1) = (-1)^e$ holds. Poisson summation now gives
$$
\theta\left(\chi,{-1 \over \tau}\right) = W_D(\chi) \left( {\tau \over i} \right)^{2e+1}
\theta(\overline\chi,\tau),
$$
with $\theta$ the series
$$
\theta(\chi,\tau) = \sum_{n \in \Z} n^e \chi(n) \exp\left({i\pi n^2\tau \over p}\right)
$$
introduced by Shimura, see e.g.\ [\Coh, Cor.\ 10.2.12]. The relation above
is most useful for $\tau = it$ with $\overline\chi(it) \not = 0$. For us,
the choice~$t=1$ has always worked, and this leads to the formula
$$
W_D(\chi) = {\theta(\chi,i)\over \theta(\overline\chi,i)}. \eqno(3.12)
$$
The key point is that we don't need that many terms of $\theta(\chi,i)$ to
use~(3.12). Indeed, the theta-series decays very rapidly, and although a
rigorous analysis appears to be complicated, we have found that only
including the first $\lceil \sqrt{p \log p} \rceil$ terms of the
theta-series works well to get an approximation to~$g(\chi)$. 

We remark that only a `rough' approximation to $g(\chi)$ is needed. 
Indeed, the $n$th power of $g(\chi)$ is known by the theorem of
Eisenstein-Weil, so that we only need to approximate
$g(\chi)$ with error in the argument less than~$\pi/n$.

\noindent
\head 4. The case $n=6$
\endhead
We now restrict to $n=6$, and fix $K = \Q(\zeta_6)$. The extension $K/\Q$ has
degree two, and since $K$ is norm-Euclidean, we can take $S = \{2,3,\infty\}$. The
localization $K_2 = \Q_2(\zeta_6)$ is the unique {\it unramified\/} degree
two extension of~$\Q_2$, and the localization $K_3 = \Q_3(\zeta_6)$ is
totally ramified of degree two.

It is a standard computation to compute the local unit groups~$U_2, U_3$,
see e.g.\ [\Neu] for an algorithm. We have
$$
U_2 / U_2^{*6} \cong \Z/6\Z \times \Z/2\Z \times \Z/6\Z \times \Z/2\Z
$$
and for convenience, we take the same generators $\alpha_1 = \pi_2 = 2, 
\alpha_2 = 3+2\zeta_6, \alpha_3 = 5+3\zeta_6, \alpha_4 = 1+2\zeta_6$ as in~[\Well].
For the other localization, we have
$$
U_3 / U_3^{*6} \cong \Z/6\Z \times \Z/3\Z \times \Z/6\Z \times \Z/3\Z
$$
with generators $\beta_1 = \pi_3 = 2\zeta_6-1=\sqrt{-3}, \beta_2 =  (1+\pi_3)^2, 
\beta_3 = 2, \beta_4 = 1+\pi_3^3$.

\subhead 4.1 Hilbert symbol
\endsubhead
\noindent
In this subsection we detail the computation of the Hilbert symbol $(x,y)_S$
on~$K_S = K_2 \times K_3 \times \C$. Firstly, the symbol is trivial on~$\C$,
so we restrict ourselves to the non-archimedean
case. The basic idea in computing $(x,y)_v$ is to write
$$
(x,y)_v = \prod_{w \not = v} (x,y)_w^{-1}
$$
by the product formula. Now, for $v \not \in 
S\cup\{ \pi_w \mid x \}$, the Hilbert symbol basically equals the power
residue symbol which we can compute by Euler's criterion,
see [\CF, Exercise 1]. The trouble lies in computation of $(x,y)_w$ for
$w \in S$ and for $w \mid x$. We will follow an idea from [\EDiss] to make
those remaining cases easy to compute as well.

We fix a place~$v\in \{2,3\}$, and let $x_1,\ldots,x_4$ be a basis
for $U_v/U_v^{*6}$. We let $w\not =v$ be the other divisor of~6. We
claim that we may assume that $x_i$ is integral and that
$$
x_i \equiv 1 \bmod K_w^{*6} 
$$
holds. To see this, we write $x_i = y_{i,w} a_{i,w}^6$ with $y_{i,w} \in \O$.
Furthermore, we let $m_w \in \Z_{>0}$ be such that $U_w^{m_w} \subset K_w^{*6}$
holds. Using the Chinese Remainder Theorem, we
choose $z_{i,w}\in\O$ with $z_{i,w} \equiv 
y_{i,w}^{-1} \bmod P_w^{m_w}$ and with $z_{i,w} \equiv 1 \bmod P_v^{m_v}$,
where $P_w$ denotes the $\O$-ideal corresponding to the valuation~$w$. In
the formula, the inverse of $y_{i,w}$ is taken in the group $K_w^*/K_w^{*6}$.
The element
$$
z_{i,w} x_i
$$
now has the desired property. We note that multiplication by $z_{i,w}$ has
not changed $x_i \bmod K_v^{*6}$.

By construction, we have~$(x_i,x_j)_w = 1$. For a place $s \not\in S$, we
have
$$
(x_i,x_j)_s = {c \legendre s} \qquad \hbox{for} \qquad c = (-1)^{s(x_i) s(x_j)}
x_i^{s(x_j)} x_j^{s(x_i)}\eqno(4.1)
$$
by [\CF, Exercise 2]. In particular, if $s(x_i) = s(x_j) = 0$, then we have
$(x_i,x_j)_s = 1$. Summarizing, we have
$$
(x_i,x_j)_v = \prod_{s} (x_i,x_j)_s^{-1} \eqno(4.2)
$$
where the product is over those places $s \not\in S$ with
either $s(x_i) \not=0$ or $s(x_j) \not=0$.  The symbols in (4.2) are easily
computed using (4.1) and the generalized Euler criterion: ${c \legendre s}$
is the unique $6$-th root of unity with
$$
{c\legendre s} \equiv c^{{N_{K/\Q}(P_s)-1\over 6}} \bmod P_s.
$$

For our choice of basis, we get $(\alpha_i,\alpha_j)_2 = \zeta_6^{a_{ij}}$ with
the matrix $A = (a_{ij})$ given by
$$
A = \pmatrix
0&3&4&0\\ 3&3&3&0\\ 2&3&0&3\\ 0&0&3&3\\
\endpmatrix
.
$$
For the localization at 3, we get $(\beta_i,\beta_j)_3 = \zeta_6^{b_{ij}}$ with
$B=(b_{ij})$ given by
$$
B = \pmatrix
3&0&3&2\\ 0&0&4&0\\ 3&2&0&0\\ 4&0&0&0\\
\endpmatrix
.
$$
The matrices satisfy $(A + A^T)\bmod 6 = (B+B^T)\bmod 6 =
0_6$. This property, which follows from $(x,y)_v (y,x)_v = 1$, is a good
check on the computation.

\subhead 4.3 Gauss sums
\endsubhead
\noindent
In this subsection we give the details on the computation 
of $g_6(1,\varepsilon, \pi)$ for a prime $\pi\in\O_S$. Assume 
first that $p=\pi \in \O_S$ is an inert prime. In this case, we
have\par
\centerline{\vbox{
\halign{$#$ \thinspace  &=\thinspace ${\displaystyle #}$\hfill\cr
g_6(1,\varepsilon,\pi) & \sum_{x,y \bmod p\Z} \varepsilon \left( {x+\zeta_6y
\over p}\right)_S e \left({x+\zeta_6y \over p}\right)\cr
& \sum_{x \not = 0} \varepsilon \left( {x+\zeta_6y\over p}\right)_S e\left(
{x+\zeta_6y\over p}\right)_S + \sum_{y \bmod p\Z} \varepsilon\left(
{\zeta_6y\over p}\right)_S e\left(\zeta_6y\over p\right).\cr
& \sum_{y \bmod p\Z} \varepsilon\left( {1+\zeta_6 y\over p} \right)_S
\sum_{x \not = 0 \bmod p} e\left( {x+\zeta_6 xy\over p}\right) - 
\varepsilon\left(\zeta_6\over p\right)_S,\cr}
}}
\noindent
where we have made the substitution $y \rightarrow xy$ and used the 
equality $\left({y\over p}\right)_S=1$ in the last line. 
We compute $\varepsilon\left({x+\zeta_6xy\over p}\right) = 
\exp(2\pi i {x(2+y)\over p})$, and derive that\par
\centerline{\vbox{
\halign{$#$ \thinspace  &=\thinspace ${\displaystyle #}$\hfill\cr
g_6(1,\varepsilon,\pi) & -\varepsilon\left({\zeta_6\over p}\right)_S +
(p-1) \varepsilon\left({1-2\zeta_6\over p}\right)_S + \sum_{y\not = -2} 
\varepsilon\left({1-\zeta_6y\over p}\right)_S\cr
& p \varepsilon\left( {1-2\zeta_6 \over p}\right)_S.\cr}
}}
\noindent
Here, the sum over all $y\not = -2$ is computed by expanding the sum in the
equality $0 = {\displaystyle \sum_{x,y\bmod p\Z}} \varepsilon\left( {1-2\zeta_6y\over p}
\right)_S$ and rearranging terms. We conclude that
$$
g_6(1,\varepsilon,\pi) = p \varepsilon\left( {1-2\zeta_6 \over p}\right)_S
= \cases p & \hbox{for\ } p \equiv 3 \bmod 4\\
          -p & \hbox{for\ } p \equiv 1 \bmod 4
\endcases
$$
holds for inert primes.

For a split prime~$\pi$ of norm~$p$, we have\par
\centerline{\vbox{
\halign{$#$ \thinspace  &=\thinspace ${\displaystyle #}$\hfill\cr
g_6(1,\varepsilon,\pi) & \sum_{x=0}^{p-1} \varepsilon\left( {x\over\pi}
\right)_S \exp(2\pi i {x\over\pi} \Tr(\pi))\cr
& \varepsilon\left( {-\overline\pi\over\pi}\right)_S^{-1} g\left(\varepsilon\left(
\cdot\over\pi\right)_S\right),\cr}
}}
\noindent
with $g$ as in~(3.11). The methode from subsection~3.4 can be used to
approximate~$g$. In this case, the Eisenstein-Weil theorem tells us
that
$$
g^3 = \pi^2\sqrt{p}\cdot\cases 1 & \hbox{if\ } p\equiv 1\bmod 4\\
                                  3 & \hbox{if\ } p\equiv 3\bmod 4
                            \endcases\eqno(4.3)
$$
holds, and we only need to approximate~$g$ with enough accuracy to
select the right cubic root of the right hand side of~(4.3).

To compute $g_6(\pi^{k-1},\varepsilon,\pi^{k})$, we need to know $g_6(1,
\varepsilon^k,\pi)$ for $k=2,3,4,5$. The value for these quantities follows
directly from the definition. Indeed, we have\par
\smallskip
\centerline{\vbox{
\halign{$#$\thinspace &=\thinspace $#$\hfill\cr
g_6(1,\varepsilon^5,\pi) & \overline{g_6(1,\varepsilon,\pi)} \varepsilon\left( 
                         {-1 \over\pi}\right)_S\cr
g_6(1,\varepsilon^4,\pi) & \overline{g_6(1,\varepsilon^2,\pi)} \varepsilon\left(
                         {-1 \over\pi}\right)_S\cr
g_6(1,\varepsilon^3,\pi) & g_2(1,\varepsilon^3,\pi)\cr
g_6(1,\varepsilon^2,\pi) & g_3(1,\varepsilon^2,\pi)\cr}}}
\smallskip\noindent
with $g_2,g_3$ the quadratic and cubic Gauss sum respectively. For the
quadratic Gauss sum we have
$$
g_2(1,\varepsilon^3,\pi) = \cases \varepsilon^3\left( {\overline{\pi}\over\pi}\right)_S \sqrt{N_{K/\Q}(\pi)} & \hbox{if\ } |\pi| \equiv 1\bmod 4\\
\varepsilon^3\left( {-\overline{\pi}\over\pi}\right)_S i\sqrt{N_{K/\Q}(\pi)} & \hbox{if\ } |\pi| \equiv 3\bmod 4,
\endcases
$$
and for the cubic Gauss sum we have
$$
g_3(1,\varepsilon^2,\pi) = g\left(\varepsilon^2\left({\cdot\over\pi}\right)_S\right).
$$
This last Gauss sum can be computed as in subsection 3.4 using the relation
$$
g^3 = N_{K/\Q}(\pi) \pi \qquad \hbox{for}\qquad \pi \equiv -1 \bmod 3.
$$

\subhead 4.4. Hypergeometric function
\endsubhead\par
\noindent
We recall that we need to evaluate the functions
$$
F_{1,2}(x) = {1\over2\pi i} \int_{\sigma} \left( {\Gamma(6s)\over\Gamma(s)}
\right) f(\pm x) x^{-s} \hbox{d}s \eqno(4.4)
$$
many times. However, we have to evaluate the function $F_2$ roughly 
$6^3=216$ times more often than~$F_1$. Hence, we will pick the function~$f$ so that
$F_2$ is especially easy to evaluate. 

Following [\Well], we propose to take $f(s) = \Gamma(1+s)^{-1}$. It can be
easily checked that
$$
F_2(x) = {1\over 6\pi} \exp\left(-\sqrt{3} {x^{1/6}\over2}\right) \sin\left({x^{1/6}\over2}\right)
$$
holds. We believe that the fact that $F_2$ is numerically easy to evaluate
makes up for the fact that (4.4) converges slower for this~$f$ than 
for~$f=1$.

The `price we pay' for the easy formula for~$F_2$ is that $F_1$ is harder
to compute. Using the residue theorem, one computes that
$$
F_1(x) = \sum_{k=1}^5 {-z^k\over \Gamma(-k/6)^2} \left({1\over k\cdot k!} 
+ \sum_{i=1}^{\infty} \left( {z^6\over 36}\right)^i { \prod_{l=1}^{i-1}(k+6l)^2\over
(k+6i-1)!}\right) \quad \hbox{with} \quad z = -x^{1/6}
$$
holds. Some remarks about this formula are in order. Firstly, although 
this series expansion for~$F_1$ converges best for small $|x|$, we have
found that even for moderately large $x \approx 10^6$ it is an efficient
way to compute~$F_1(x)$. However, care must be taken to perform all 
computations with high precision. For $x \approx 10^6$, a precision of 200
bits sufficed for us. Secondly, it is best to make a table of the
quotients
$$
{\prod_{l=1}^{i-1} (k+6l)^2\over (k+6i-1)!}
$$
instead of computing them on the fly. The first few hundred values for~$i$
suffice. Finally, we have found that the series expansion $\sum_{i=1}^{\infty}$
that occurs 
converges smoothly, so that simply checking that the summand is less than
some chosen bound suffices to approximate this series by a partial sum.\par
\ \par
\subhead 4.5 Transition matrix
\endsubhead

We use formula (3.10) to compute the coefficients~$c_{jw}$. This is relatively
straightforward, albeit it technical. Since mistakes are easy to make in
this part, we give some details on the computation and on checks one can
do to make sure the matrix~$T$ is correct. 

The local differents are $d_3 = 1$ and $d_2 = 0$. 
The local Gamma functions $\Gamma_v$ have coefficients
in the ring~$\Z[\zeta_{72}]$. It turns out that the coefficients~$T_{ij}$
themselves have coefficients in~$\Z[\zeta_{36}]$ already. Since computing
(3.10) involves taking various quotients in $\Z[\zeta_{36}]$, it is 
important for the practical performance to view $T_{ij}(r,s) \in 
\Z[\zeta_{36}]((1/3)^s,(1/4)^s)$.

Putting $V = (1/3)^s, W = (1/4)^s$, we have found that $T_{ij}$ always has 
denominator~$V^5 W^3$. The coefficients $c_{jw}= c_{j,(w_1,w_2)}$ are only
non-zero for $(w_1,w_2) \in \{-5,\ldots,4\} \times \{-3,\ldots,5\}$.
An excellent check for our computations is to compute
the product
$$
T_{ij}(r,s) T_{ij}({r,-s}).
$$
By the functional equation~(3.1), this product is a {\it diagonal matrix\/}.
Simply by computing a few coefficients of the product we can check if the
matrix has been computed correctly. 

\subhead 4.6 Computing the integrals
\endsubhead\par
\noindent
To compute the sums (3.7) and (3.9), we loop over all $\O_S$-ideals $I$
and for each $I$, we compute a generator $c$ that is 
contained in~$V$. Using the formulas from subsections 3.1 and 3.2 we
then compute the contribution from~$c$ to the respective sums. However,
since the Gauss sum $g_6(1,\varepsilon,c)$ is {\it zero\/} for $c$ that are
not squarefree, we can restrict our attention to squarefree ideals~$I$ and
the~$I$ that have a non-trivial gcd with $(r)$.

To simplify the exposition, we restrict to the case that $r$ is 
coprime to~6. We then have the following algorithm
for approximating the Fourier coefficient $\tau(r,V) = N_{K/\Q}(r)^{1/12}\text{\rm Res}_{s=1/n}\sum_{\eta\in V}\Psi(r,s,\eta).$

\proclaim {4.1 Algorithm}
\endproclaim
\noindent
{\bf Input.} An element $r \in \O$ that is coprime to 6, 
a control parameter $x>0$, a bound $B>0$,
a choice of representatives~$V$, and a precision bound~$X$.
\ \par\noindent
{\bf Output.} An approximation $(3.7)-(3.9)$ to $\tau(r,V)$ coming from taking all ideals of norm up to $B$ into account, and by performing all computations
with precision~$X$.\par
\ \par
\parindent=4pc
\item{\sl Step 1.\/} Fix the embedding $\varepsilon(\zeta_6) = \exp(2\pi i/6)$.Compute and store~$g_6(1,\varepsilon,\pi)$ for all prime ideals
$(\pi)\subset\O$ with $7 \leq |\pi| \leq B$ using the method from subsection~4.3.
\item{\sl Step 2.\/} Compute and store the coefficients~$\varepsilon(c_{j,w})$ 
of the
transition matrix using formula (3.10) for all $1 \leq j
\leq 216$ and $w \in \Z^2 \cap [-5,4] \times [-3,5]$. 
\item{\sl Step 3.\/} Compute and store the coefficients $a(m), b(m)$ for
all $m \geq 1$ until both $a(m), b(m) \leq X$ using formula (3.6).
\item{\sl Step 4.\/} Set $V_1 \leftarrow 0$, $V_2 \leftarrow 0$. Initialize
an empty list $L \subset \O_S \times \C \times \Z_{>1}$. (We will add
triples $(\alpha,g_6(r,\varepsilon,\alpha),N_{K/\Q}(\alpha))$ to $L$ later.)
Initialize an empty list $M \subset \Z_{>1} \times \C$. (We will add 
values $F_1(\cdot)$ to~$M$ later for all norms we encouter.)
\item{\sl Step 5.\/}  ({\sl Constant term\/}) Determine $j \in  \{1,\ldots\,\!,216\}$ with $\eta_j = 1
\in K_S^*/K_S^{*6}\O_S^*$. For all $w = (w_1,w_2) \in
\{-5,\ldots\,\!, 4\}\times\{-3,\ldots\,\!,5\}$ do the following.
\itemitem{\sl (a)\/} Set $f_{2,w} \leftarrow 0$.
\itemitem{\sl (b)\/} For all $k \geq 1$ do the following.
\itemitem{\sl (c)\/} Compute $s = F_2(y_1 x^{-1} (1/4)^{w_1} (1/3)^{w_2} k^6) b(k^6)$. If $|s| < X$, goto step 5d, else set $f_{2,w} \leftarrow f_{2,w} +s$
and repeat.
\itemitem{\sl (d)\/} Set $V_2 \leftarrow V_2 + f_{2,w}\cdot c_{j,(w_1,w_2)}.$
\item{\sl Step 6.\/}  ({\sl Constant term\/}) Set $f_1 \leftarrow 0$. For all $k \geq 1$ do 
                     the following.
\itemitem{\sl (a)\/} Compute $s = F_1(y_1 x k^6) a(k^6)$. If $|s| < X$, goto
step 6b, else set $f_{1} \leftarrow f_{1}+s$ and repeat.
\itemitem{\sl (b)\/} Set $V_1 \leftarrow V_1 + f_1$. Add $(1,f_1)$ to $M$.
\item{\sl Step 7.\/}  ({\sl Constant term\/}) Add $(1,1,1)$ to~$L$.
\item{\sl Step 8.\/} For all primes ideals $(\pi)\subset\O_S$ with $|\pi| < B$ (ordered 
by norm) do the following.
\item{\sl Step 9.\/} For all $(\alpha,g_6(r,\varepsilon,
\alpha),N_{K/\Q}(\alpha)) \in L$ do the following. 
\item{\sl Step 10.\/} Find $k\in \{1,\ldots\,\!6\},j \in \{1,\ldots\,\!,216\}$ 
with $\zeta_6^k a \pi \sim \eta_j \in V$. Set $a \leftarrow \zeta_6^k a$. Set $N \leftarrow N_{K/\Q}(a\pi)$.
\item{\sl Step 11.\/} Compute $g_6(r,\varepsilon,a\pi)$ using Lemma 3.2. 
\item{\sl Step 12.\/} If $N \le B/p$, add $(a\pi,
                      g_6(r,\varepsilon,a\pi),N)$ to~$L$.
\item{\sl Step 13.\/} For all $w = (w_1,w_2) \in \{-5,\ldots\,\!, 4\}
                     \times\{-3,\ldots\,\!,5\}$ do the following.
\itemitem{\sl (a)\/} Set $f_{2,w} \leftarrow 0$.
\itemitem{\sl (b)\/} For all $k \geq 1$ do the following.
\itemitem{\sl (c)\/} Compute $s = F_2(y_1 x^{-1} (1/4)^{w_1} (1/3)^{w_2} k^6 N) b(k^6)$. If $|s| < X$ goto step 13d, else set $f_{2,w} \leftarrow f_{2,w} +s$
and repeat.
\itemitem{\sl (d)\/} Set $V_2 \leftarrow V_2 + f_{2,w}\cdot c_{j,(w_1,w_2)} 
                     {g_6(r,\varepsilon,a\pi) \over N}.$
\item{\sl Step 14.\/} If $(N,x)$ is present in~$M$, set $f_1 \leftarrow x$.
Else, for all $k \geq 1$ do the following.
\itemitem{\sl (a)\/} Compute $s = F_1(y_1 x N k^6) a(k^6)$. If $|s| < X$, goto
step 14b, else set $f_{1} \leftarrow f_{1}+s$ and repeat.
\itemitem{\sl (b)\/} Add~$(N,f_1)$ to $M$.
\item{\sl Step 15.} Set $V_1 \leftarrow V_1 + f_1 {g_6(r,\varepsilon,a\pi)
\over N}$. Go to step 9.
\item{\sl Step 16.} If $|\pi| < B$, go to Step 8.
\item{\sl Step 17.} Set $\tau \leftarrow (V_1-6\sqrt{3}V_2) x^{1/6}$. Set
$\tau \leftarrow \tau\cdot \Gamma(1/6)\Gamma(7/6) y_1^{1/6}N_{K/\Q}(r)^{1/12}$. Return~$\tau$.
\par\ \par\noindent
\parindent=1pc
We make some remarks about the algorithm. Firstly, Step~1 of the algorithm
is independent of~$r$ and one should store the Gauss sums in a file once
and for all if we are computing several coefficients~$\tau(r,V)$. Secondly,
the coefficients $c_{j,w}$ in Step 2 only depend on $r \in K_S^*/K_S^{*6}$.
Although this quotient group has size $6^6 = 46656$, we still store the
matrix in a file. This is particularly convenient if we are running the
algorithm for various choices of~$x$, $X$ and~$L$. 

Since the values $f_1,f_2$ we compute get multiplied by a Gauss sum later
on, we only need to consider squarefree elements of~$\O_S$. Our loop over
all squarefree elements of norm at most~$B$ is basically a variant of the
sieve of Eratosthenes. The purpose of the list~$M$ is simply 
to avoid some (costly) evaluations of the function~$F_1$.

Finally, we remark that although all computations are done with 
precision~$X$ in the algorithm, this does not mean that the {\it output\/}
is correct with the same precision. Indeed, the precision of the output
depends on the combined choice of $B$ and $x$ as well. We refer to [\EDiss]
for an analysis of the convergence properties of sums~(3.7) and~(3.9).

\noindent
\head 5. The conjecture
\endhead
\noindent
The first thing we need to decide for actual computations is which set
of representatives~$V$ to use. After trying several possibilities, we
have found that the set
$$
V = \left\{ \eta\in K_S^*/K_S^{*6} \mid \alpha_1 = \beta_1 = \beta_2 = 0 
\hbox{\ and\ } \cases \alpha_3 = 0, & \hbox{if\ } \alpha_2 = 0\\
                      \alpha_3 = 1, & \hbox{if\ } \alpha_2 = 1
                \endcases
    \right\}
$$
gives the cleanest results. We note that this is the set $V_2$ considered
by Wellhausen in his thesis~[\Well]. The set~$V$ can be characterized as
follows:
$$
\{ x \in \O \mid \gcd(x,6) = 1, x \sim \eta \in V \} = \{ x \in \O \mid 
x \equiv y \bmod 12\}
$$
with $y\in$
$$
\{ 1, 5, 4+3\zeta_6, 8+3\zeta_6,1+6\zeta_6,5+6\zeta_6,1+9\zeta_6,2+9\zeta_6,5+9\zeta_6, 7+9\zeta_6,10+9\zeta_6,11+9\zeta_6 \}. \eqno(5.1)
$$
The choice of~$V$ is motivated by the following property:
$$
\forall v,w \in V: (v,w)_6 = 1 \qquad \hbox{or} \qquad (v,-w)_6 = 1,
$$
which can be proved easily. This property does not hold for the perhaps 
easier looking choice $V'$ characterized by
$$
\{ x \in \O \mid \gcd(x,6) = 1, x \sim \eta \in V' \} = \{ x \in \O \mid 
x \equiv 1 \bmod 3\}.
$$

The fact the Hilbert symbol is particularly easy on~$V$ has consequences
for the coefficient $\tau(\pi^4,V)$. The theory of {\it Hecke operators\/} is
used in the proof of the following lemma.

\proclaim{5.1 Lemma} Let $\pi\in\O_S$ be prime with $|\pi| \equiv 1 \bmod 4$
and $\pi \equiv y \bmod 12$ for some $y$ in set (5.1),
and let $V$ be as above.
Then we have
$$
{\tau(\pi^4,V)\over \tau(1,V)} = |\pi|^{-1/2} \overline{g_6(1,\varepsilon,\pi)}.
$$
\endproclaim\noindent
{\bf Proof.} The following holds for $\rho(\pi^4,\eta) = 
                    \hbox{Res}_{s=1/6} \psi(\pi^4, s,\eta)$:
$$
\rho(\pi^4,\eta) = |\pi|^{-5/6} g_6(1,\varepsilon^5,\pi) 
                    \varepsilon((-\eta,\pi^5)_S) \rho(1,\eta \pi^{-5})
$$
see~[\KP]. We sum over all $\eta\in V$ and renormalize to $\tau$ to obtain
$$
\tau(\pi^4,V) = |\pi|^{-1/2} g_6(1,\varepsilon^5,\pi)
     \sum_{\eta = \alpha\pi^5\in V} \varepsilon((\alpha,\pi^5)_S)\rho(1,\alpha)
$$
where we have used the equality $(-\pi^5,\pi^5)_6 = 1$ in the
sum. By replacing $\pi$ by $-\pi$ if necessary, we may assume that
$\alpha \in V$. The lemma follows from checking that for $|\pi| \equiv 1 \bmod 4$
and $\alpha\pi\in V$, the equalities $(\alpha,\pi^5)_6$ hold for 
all $\eta\in V$.  \hfill\endproof\par
\ \par\noindent
We remark that the proof hinges on the special property of~$V$. If we 
change~$V$, then the lemma need not be true. Furthermore, for $|\pi| \equiv 3
\bmod 4$, not all the Hilbert symbols $(\alpha,\pi^5)_S$ are trivial. The
lemma is false in this case as well. 

To see what happens for $|\pi| \equiv 3 \bmod 4$, we use Algorithm 4.1. We
have found that the choices
$$
B = 10^8 \qquad\qquad x = 1/300 \qquad\qquad X = 10^{-20}
$$
work very well. We compute and store all Gauss sum for primes up to 
norm~$10^8$. This computation is highly parallelizable, and it is of great
help here to have a cluster of CPU's available. We used the method from
subsection~3.4 to compute the individual Gauss sums, noting that the
equality
$$
g_6(1,\varepsilon,\overline\pi) = (-1,\pi)_6 g_6(1,\varepsilon,\pi)
$$
saves us half the computations. 

Lemma 5.1 is a very good test for the implementation, since
a small mistake in the implementation will cause the equality in Lemma 5.1
to be false. Furthermore, the output of the algorithm should be roughly
independent of~$x$. By letting $x$ vary over $1/500,1/400,1/300,1/200,1/100$,
we can check that the algorithm is performing correctly. The first quantity
to compute is $\tau(1,V)$. In agreement with [\Well], we find that
$$
\tau(1,V) \approx 0.1358547858696091.
$$
By letting~$x$ vary and checking the independence of~$x$ in the computations,
we are confident that the expression above is correct up to 16 decimal digits. We
remark that for other choices of~$V$, the `constant term' $\tau(1,V)$ need
not be real.

\proclaim{5.2. Conjecture} Let $\pi\in\O_S$ be prime with $|\pi| \equiv 3
\bmod 4$
and $\pi \equiv y \bmod 12$ for some $y$ in set (5.1). Then we have
$$
{\tau(\pi^4,V)\over \tau(1,V)} = |\pi|^{-1/2} {(-1,\pi)_6 \overline{g_6(1,\varepsilon,\pi)}\over \sqrt{3}}.
$$
\endproclaim\noindent
{\bf Evidence.} This conjecture is purely based on computational evidence. 
In fact, since the norm of $\pi^4$ grows rather quickly, we only computed
a few cases. The conjecture is correct for $|\pi| = 7,19,31$ for several
decimal digits. To {\it prove\/} this conjecture, one should examine the
relations between the $\rho(\pi^4,\eta)$ for varying $\eta$ more closely.
\hfill\endproof\par
\ \par\noindent
We remark that for this choice of~$V$, the coefficients $\tau(\pi^4,V)$ are
in agreement with the general philosophy explained in the introduction. We
now move on to the coefficients $\tau(\pi^2,V)$. In this case, the Hecke
operators relate $\rho(\pi^2,\eta)$ to itself. However,
we can still derive the following.

\proclaim{5.3. Lemma} Let $\pi\in\O_S$ be prime with $|\pi| \equiv 1 \bmod 12$, and $\pi \equiv y \bmod 12$ for some $y$ in set (5.1). If 
we have $\tau(\pi^2,V) \not = 0$, then ${\overline\pi\legendre\pi}_6 = 1$.
\endproclaim
\noindent
{\bf Proof.} The Hecke operators now give
$$
\rho(\pi^2,\eta) = |\pi|^{-5/6} g_2(1,\varepsilon^3,\pi) \varepsilon(-\eta,
\pi^3)_S \rho(\pi^2,\eta\pi^{-3})
$$
see [\KP]. Analagous to the proof of Lemma~5.1 we derive that
$$
\tau(\pi^2,V) = |\pi|^{-1/2} g_2(1,\varepsilon^3,\pi) \tau(\pi^2,V)
$$
holds for $|\pi| \equiv 1 \bmod 12$. Furthermore, we 
have $g_2(1,\varepsilon^3,\pi) = {\overline\pi\legendre\pi}_6 \sqrt{N_{K/\Q}(\pi)}$ in
this case. The lemma follows.\hfill\endproof\par
\ \par\noindent
We caution that the converse of the lemma does not hold. In our computations
we have found several cases where $\tau(\pi^2,V) \approx 0$ even though
${\overline\pi\legendre\pi}_6 = 1$. Specifically, we conjecture that
$\tau(\pi,V) = 0$ for 
$$
|\pi| = 37, 313, 373, 661, 769.
$$
These four norms are the only norms less than 1300 for which ${\overline\pi
\legendre\pi}_6 = 1$ and $\tau(\pi^2,V) \approx 0$. We have not been able
to determine a pattern in this small set of primes.

\proclaim{5.4. Conjecture} Let $\pi\in\O_S$ be prime with $|\pi| \equiv 1
\bmod 12$, and $\pi \equiv y \bmod 12$ for some $y$ in set (5.1). If $\tau(\pi^2,V) \not = 0$, then we have
$$
{\tau(\pi^2,V) \over \tau(1,V)} = \zeta_6^{g(\pi)}
          {2 g_3(1,\varepsilon^2,\overline\pi) \over |\pi|^{1/2}}
$$
for some value $g(\pi) \in \{0,\ldots\,\!,6\}$ satisfying $g(\pi)+
g(\overline\pi) \equiv 0 \bmod 6$.
\endproclaim\noindent
{\bf Evidence.} The support of this conjecture is numerical. We have
approximated $\tau(\pi^2,V)$ for all $|\pi| < 1300$. We list the
values for $g(\pi)$ below for all $|\pi| < 10^3$.\par
\smallskip
\centerline{\vbox{\offinterlineskip 
\halign{\hfil$#$\ \hfil&\vrule #&
\ \hfil$#$\ \hfil&\vrule#&\strut\hfil \ $#$\hfil&\qquad#& \hfil$#$\ &\vrule#&\ \hfil$#$\ \hfil&\vrule#&\hfil\ $#$\hfil &\qquad#& \hfil$#$\ \hfil&\vrule#&\ \hfil$#$\ \hfil&\vrule#&\hfil\  $#$\hfil\cr
N(\pi)&&\pi&&g(\pi)&&N(\pi)&&\pi&&g(\pi)&&N(\pi)&&\pi&&g(\pi)\cr
&height2pt&&height2pt&\omit&&&height2pt&&height2pt&&&&height2pt&&height2pt\cr
\multispan5{\hrulefill}&&\multispan5{\hrulefill}&&\multispan5{\hrulefill}\cr
&height2pt&&height2pt&\omit&&&height2pt&&height2pt&&&&height2pt&&height2pt\cr
61&&-9\zeta_6+4&&5&&433&&-24\zeta_6+13&&0&&853&&-27\zeta_6+31&&1\cr
157&&-12\zeta_6+13&&2&&577&&27\zeta_6-8&&5&&877&&3\zeta_6+28&&4\cr
193&&-9\zeta_6+16&&1&&601&&-24\zeta_6+25&&3&&937&&3\zeta_6-32&&2\cr
349&&-3\zeta_6-17&&1&&613&&9\zeta_6+19&&2&&977&&36\zeta_6-23&&3\cr
397&&12\zeta_6-23&&1&&673&&21\zeta_6-29&&4&&&\cr}
}}
\smallskip\noindent
We have not been able to find a pattern in the exponents~$g(\pi)$.\hfill\endproof\par
\ \par\noindent
For the $\pi^2$-case, it remains to consider the inert primes and the
primes of norm congruent to $7 \bmod 12$. For an inert prime~$\pi$, we 
remark that although $|\pi| \equiv 1 \bmod 4$, the residue symbol in Lemma~5.3
is undefined and the proof therefore does not follow through. We have
the following conjecture.

\proclaim{5.5. Conjecture} Let $\pi\in\O_S$ be an inert prime with 
$\tau(\pi^2,V) \not = 0$, and $\pi \equiv y \bmod 12$ for some $y$ in set (5.1). Then we have
$$
{\tau(\pi^2,V) \over \tau(1,V)} = -{2 \sqrt{|\pi|}\over |\pi|^{1/2}}.
$$
\endproclaim\par\noindent
{\bf Evidence.} We have computed the coefficients~$\tau(\pi^2,V)$ for 
$\pi = 5,11,\ldots,89$. We have $\tau(\pi^2,V) \approx 0$ for
$$
\pi = 5,17,29,41,53,59,89
$$
and the conjecture is true for the other cases $\pi = 11,23,47,71,83$ with
several decimal digits precision. From this data one can furthermore
conjecture that $\tau(\pi^2,V) = 0$ for $|\pi| \equiv 5 \bmod 12$. Note
that the prime 59 contradicts the converse statement.\hfill\endproof\par
\ \par\noindent
We remark that the $\sqrt{|\pi|}$ in Conjecture 5.5 equals
the cubic Gauss sum for $\pi$, just like in Conjecture 5.4.

\proclaim{5.6. Conjecture} Let $\pi\in\O_S$ be prime with $|\pi| \equiv
7 \bmod 12$, and $\pi \equiv y \bmod 12$ for some $y$ in set (5.1). Then we have
$$
{\tau(\pi^2,V) \over \tau(1,V)} = \zeta_{12}^{h(\pi)} {2 g_3(1,\varepsilon^2,
\overline\pi)\over \sqrt{3} |\pi|^{1/2}}
$$
with $h(\pi) \in \{1,3,5,7,9,11\}$ satisfying $h(\pi)+h(\overline\pi) \equiv
0 \bmod 12$.
\endproclaim
\noindent
{\bf Evidence.} We have approximated the coefficients $\tau(\pi^2,V)$ for
all $|\pi| < 8000$. We list the first few values for $h(\pi)$ below.\par
\smallskip
\centerline{\vbox{\offinterlineskip 
\halign{\hfil$#$\ \hfil&\vrule #&
\ \hfil$#$\ \hfil&\vrule#&\strut\hfil \ $#$\hfil&\qquad#& \hfil$#$\ \hfil&\vrule#&\ \hfil$#$\ \hfil&\vrule#&\hfil\ $#$\hfil &\qquad#& \hfil$#$\ \hfil&\vrule#&\ \hfil$#$\ \hfil&\vrule#&\hfil\  $#$\hfil\cr
N(\pi)&&\pi&&h(\pi)&&N(\pi)&&\pi&&h(\pi)&&N(\pi)&&\pi&&h(\pi)\cr
&height2pt&&height2pt&\omit&&&height2pt&&height2pt&&&&height2pt&&height2pt\cr
\multispan5{\hrulefill}&&\multispan5{\hrulefill}&&\multispan5{\hrulefill}\cr
&height2pt&&height2pt&\omit&&&height2pt&&height2pt&&&&height2pt&&height2pt\cr
7&&3\zeta_6-2&&9&&79&&3\zeta_6-10&&7&&163&&3\zeta_6+11&&3\cr
19&&3\zeta_6+2&&11&&103&&-9\zeta_6-2&&1&&199&&15\zeta_6-13&&11\cr
31&&6\zeta_6-1&&3&&127&&6\zeta_6-13&&7&&211&&15\zeta_6-1&&7\cr
43&&-6\zeta_6+7&&1&&139&&3\zeta_6-13&&7&&223&&6\zeta_6+11&&3\cr
67&&-9\zeta_6+7&&9&&151&&-9\zeta_6+14&&7&&271&&-9\zeta_6+19&&1\cr}
}}
\smallskip\noindent
We have not been able to find a pattern in the exponents~$h(\pi)$. Since 
$h(\pi)$ appears to be always odd, we can replace the $\sqrt{3}$ in the
denominator of the conjecture by $\zeta_{12}^{11}(1+\zeta_6)$ to force
$h(\pi)$ to be even and $\zeta_{12}^{h(\pi)}$ is then a {\it sixth\/} root of
unity. However, complex conjugation does not act nicely in this case. This
is the reason we have stated the conjecture with a $\sqrt{3}$ instead
of $1+\zeta_6$. \hfill\endproof\par
\ \par
\noindent
We proceed with the investigation of $\tau(\pi,V)$. We believe that the
quantity 
$$
\left( {\tau(\pi,V) \over \tau(1,V)}\right)^2
$$ 
has interesting algebraic properties. 

\proclaim{5.7. Conjecture} Let $\pi\in\O_S$ be prime with $\pi \equiv y \bmod 12$ for some $y$ in set (5.1).  If $|\pi| \equiv 1
\bmod 4$, then we have
$$
\left( {\tau(\pi,V) \over \tau(1,V)}\right)^2 = \zeta_6^{k(\pi)} {g_3(1,\varepsilon^2,\overline\pi)\over |\pi|^{1/2}} 3^{l(\pi)} (1-3\zeta_6)^{m(\pi)} (-2+3\zeta_6)^{n(\pi)}
$$
with $k(\pi) \in \{1,\ldots,6\}$ satisfying $k(\pi) + k(\overline\pi) =  0
\bmod 6$. We have $l(\pi) \in \{-1,0\}$. The elements $1-3\zeta_6,-2+3\zeta_6$ have norm 7, and we have
$m(\pi), n(\pi) \in \{0,2\}$ with the restriction that they cannot both
be equal to 2. If $|\pi| \equiv 3 \bmod 4$, then we have
$$
\left( {\tau(\pi,V) \over \tau(1,V)}\right)^2 = \zeta_6^{k(\pi)} {g_3(1,\varepsilon^2,\overline\pi)\over |\pi|^{1/2}} 3^{l(\pi)} (1+3\zeta_6)^{m(\pi)} (4-3\zeta_6)^{n(\pi)}
$$
with the same restrictions on $k,l,m,n$. The elements $1+3\zeta_6,4-3\zeta_6$
both have norm 13.
\endproclaim
\noindent
{\bf Evidence.} In the case $|\pi| \equiv 1 \bmod 4$, we have computed 
$\tau(\pi,V)$ for all $|\pi| < 900$. The norms where an element of norm 7
appears in $\tau$ are
$$
73, 193, 241, 349, 373, 421, 613, 661,709,757,829.
$$
For these norms, we have $l(\pi) = -1$. The norms with $l(\pi) = 0$ are
$$
97, 229,313,457,577,877.
$$
The conjecture is on thinner ice for $|\pi| \equiv 1 \bmod 4$. In this case,
our implementation is not entirely independent of the parameter~$x>0$, which
has the practical impact that we can only rely on very few digits. In his
thesis, Wellhausen computed $\tau(\pi,V)$ for all $|\pi| < 100$ and the
only $|\pi|$ where the element of norm 13 appears is $79$. The norms with
$l(\pi) = 0$ are
$$
19,31.
$$
We have not found a pattern in the exponents~$k(\pi)$, nor a condition
when the elements of norm 7,13 appear. \hfill\endproof\par

\head Acknowledgement
\endhead
\noindent
It is a great pleasure to thank Professor Patterson for many helpful discussions.

\enddocument